\documentclass[11pt,a4]{article}

\usepackage[utf8]{inputenc}
\usepackage[english,main=russian]{babel}

\usepackage{amsfonts, amsmath, amsthm}
\usepackage{amssymb}
\usepackage{url}
\usepackage[
sorting=none
]{biblatex}

\usepackage{graphicx}
\usepackage{subcaption}

\usepackage{mathtools}
\usepackage{makecell}
\usepackage{tikz}
\usepackage{dsfont}
\usepackage{secdot}

\usepackage{algorithm}
\usepackage{algorithmic}

\sectiondot{section}
\sectiondot{subsection}
\theoremstyle{definition}

\newcommand{\norm}[1]{\left\lVert#1\right\rVert}

\newtheorem{definition}{Определение}
\newtheorem{remark}{Замечание}
\newtheorem{hypothesis}{Гипотеза}

\newcommand\myeq{\mathrel{\stackrel{\makebox[0pt]{\mbox{\normalfont\tiny def}}}{=}}}

\DeclareMathOperator*{\argmin}{arg\,min}

\newcounter{algorithm_heading}
\newenvironment{algorithm_heading}[1][]{\refstepcounter{algorithm_heading}\begin{center}\bf}{\end{center}}
\allowdisplaybreaks

\def\at2#1{{\color{black}#1}}

\begin{document}

\textbf{УДК} 519.85

\begin{center}
\textbf{Эвристический адаптивный быстрый градиентный метод в задачах стохастической оптимизации\footnote{Работа  А.И. Тюрина поддержана грантом РФФИ 19-31-90062 Аспиранты.}}

{\bf
\copyright\,2019 г.\,\,
А. В. Огальцов$^*$,
А. И. Тюрин$^*$
}

(*101000 Москва, ул. Мясницкая, 20, НИУ ВШЭ;)
 
e-mail: atyurin@hse.ru

Поступила в редакцию: **.**.**** г.
\end{center}

\renewcommand{\abstractname}{\vspace{-\baselineskip}}
\begin{abstract}
В данной работе приводится эвристический стохастический адаптивный ускоренный градиентный метод. Мы показываем, что на практике предложенный алгоритм имеет высокую скорость сходимости по сравнению с популярными сейчас методами оптимизации. Кроме того, мы приводим обоснование нашего метода и описываем трудности, которые не позволяют на данный момент получить оптимальные оценки для предложенного алгоритма.
\end{abstract}

\textbf{Ключевые слова:} быстрый градиентный метод, стохастическая оптимизация, адаптивная оптимизация.

\section{Введение}
В данной работе представлен эвристический алгоритм стохастической оптимизации, который является обобщением быстрого градиентного метода \cite{nesterov2010introductory}. Методы стохастической оптимизации в последнее время являются популярными по той причине, что они позволяют уменьшать сложность подсчета градиента, что является очень важным, так как существуют примеры функций, в которых невозможно за разумное время подсчитать градиент оптимизируемой функции хотя бы в одной точке \cite{goodfellow2016deep}, \cite{krizhevsky2012imagenet}. Помимо этого, в задачах стохастической оптимизации \cite{gasnikov2016netrivialnosti} в силу формулировки самой задачи единственный разумный способ получить направление спуска --- это сэмлирование векторов, которые имеют несмещенные по отношению к истинному градиенту направления. Данный подход также популярен в глубинном обучении \cite{krizhevsky2012imagenet}, в задачах, которые имеют вид суммы функций \cite{gasnikov2016netrivialnosti}. Стоит отметить, что сложность подсчета стохастического градиента, как правило, можно регулировать.

Задача по поиску адаптивного быстрого градиентного метода с стохастическим градиентом по-прежнему является актуальной. На сколько нам известно данная проблема является не решенной. Отметим, что в данном направлении активно идет работа. В работе \cite{bach2019universal} предложили метод, который одновременно хорошо работает на гладких и негладких задачах, однако, данный метод не является ускоренным, шаг метода выбирается таким образом, что он не возрастает, в то время, как в нашем методе шаг метода подбирается адаптивно, может возрастать и убывать. В работе \cite{vaswani2019painless} предложили объединить стохастический неускоренный градиентный метод и правило Армихо \cite{nocedal2005numerical}, в результате авторы доказали, что их метод имеет скорость сходимости неускоренного градиентного метода, но это им удалось, используя специальные предположения о оптимизационной задаче. 

К сожалению, нам пока так и не удалось получить теоретические оценки для предложенного метода. В данной работе приведены эксперименты (см. раздел \ref{comp_exp}), которые показывают, что в задачах машинного обучения с логистической регрессией, нейроной сетью и сверточной нейронной сетью предложенный алгоритм имеет более быструю скорость сходимости, чем Adagrad \cite{duchi2011adaptive} и Adam \cite{adam}. Кроме того, в разделе \ref{explanation} мы приводим некоторое обоснование нашего алгоритма и объясняем некоторые детали, которые нам не позволяют до конца доказать сходимость метода в оптимальной форме.

\section{Адаптивный быстрый стохастический градиентный метод} \label{sec:mmtDLST}

Опишем сначала общую постановку задачи выпуклой оптимизации \cite{nesterov2010introductory}. 
Пусть определена функция $f(x): Q \longrightarrow \mathds{R}$ и дана произвольная норма $\norm{\cdot}$ в $\mathds{R}^n$. Сопряженная норма определяется следующим образом:
\begin{gather*}\norm{\lambda}_* \myeq \max\limits_{\norm{\nu} \leq 1;\nu \in \mathds{R}^n}\langle \lambda,\nu\rangle,\,\,\,\forall \lambda \in \mathds{R}^n.\end{gather*}
Будем полагать, что
\begin{enumerate}
	\item $Q \subseteq \mathds{R}^n$, выпуклое, замкнутое.
	\item $f(x)$ -- непрерывная, гладкая и выпуклая функции на $Q$.
	\item $f(x)$ ограничена снизу на $Q$ и достигает своего минимума в некоторой точке (необязательно единственной) $x_* \in Q$.
\end{enumerate}

Рассмотрим следующую задачу оптимизации:
\begin{align*}
f(x) \rightarrow \min_{x \in Q}.
\end{align*}

Введем два понятия: прокс-функция и дивергенция Брэгмана \cite{gupta2008bregman}.
\begin{definition}
$d(x):Q \rightarrow \mathds{R}$ называется прокс-функцией, если $d(x)$ непрерывно дифференцируемая на $\textnormal{int }Q$ и $d(x)$ является 1-сильно выпуклой относительно нормы $\norm{\cdot}$ на $\textnormal{int }Q$.
\end{definition}
\begin{definition}
Дивергенцией Брэгмана называется 
\begin{align*}
V(x,y) \myeq d(x) - d(y) - \langle\nabla d(y), x - y\rangle,
\end{align*}
где $d(x)$ --- произвольная прокс-функция.
\end{definition}

\begin{definition}
\leavevmode
\label{defdeltaLstoch}
Стохастическим $(\delta, L)$-оракулом будем называть оракул, который на запрашиваемую точку $y \in Q$ дает пару $\left(f_\delta(y), \nabla f_\delta(y;\xi)\right)$ такую, что
\begin{gather*}
0 \leq f(x) - f_\delta(y) - \langle\nabla f_\delta(y), x - y\rangle \leq \frac{L}{2}\norm{x - y}^2 + \delta ,\,\,\, \forall x \in Q,
\end{gather*}
\begin{gather*}
\mathbb{E}\left[\nabla f_\delta(y;\xi)\right] = \nabla f_\delta(y),\,\,\, \forall y \in Q,
\end{gather*}
\begin{gather*}
\mathbb{E}\left[\exp\Bigg(\frac{\norm{\nabla f_\delta(y;\xi) - \nabla f_\delta(y)}^2_*}{\sigma^2}\Bigg)\right] \leq \exp(1),\,\,\, \forall y \in Q,
\end{gather*}
где $\sigma^2 > 0$ и $\xi$~--- произвольная случайная величина.
\end{definition}

Отметим, что концепция стохастичксого $(\delta, L)$--оракула основана на концепции $(\delta, L)$--оракула \cite{devolder2013first}, \cite{devolder2014first}, \cite{devolder2013exactness}.

\begin{remark}
\label{remark_stoch_fun}
В определении \ref{defdeltaLstoch} можно дополнительно потребовать случайность не только градиента, но и самой функции, то есть, чтобы стохастический $(\delta, L)$-оракул возвращал вместо неслучайного $f_\delta(y)$ значения функции некоторое случайное значение $f_\delta(y;\xi)$.
\end{remark}

Определим константу $R_Q$ такую, что
\begin{align*}
R_Q \geq \max_{x,y \in Q}\norm{x - y}.
\end{align*} Будем считать, что $R_Q < \infty$.

В методе, который мы предложим далее, будем оценивать истинный градиент на каждом шаге с помощью некоторого количества $\nabla f_\delta(y;\xi_j),\, j \in [1\dots m_{k+1}]$, используя технику mini-batch (см, например, \cite{li2014efficient}).

Обозначим
\begin{gather}
\label{stoch_gradient}
\widetilde{\nabla}^{m_{k+1}} f_\delta(y) \myeq \frac{1}{m_{k+1}}\sum_{j=1}^{m_{k+1}}\nabla f_\delta(y;\xi_j).
\end{gather}

Пусть $\kappa$ -- константа регулярности из \cite{juditsky2008large} для $(\mathds{R}^n,\norm{\cdot}_*)$. Для некоторых простых случаев верно следующее: если $\norm{\cdot}_* = \norm{\cdot}_2$, то $\kappa = 1$. Более того, если $\norm{\cdot}_* = \norm{\cdot}_q$, $q \in [2,\infty]$, то $\kappa = \min\left[q-1,2\ln(n)\right]$. Сделаем следующее обозначение: $\widetilde{\Omega} \myeq 2\kappa + 4\Omega\sqrt{\kappa} + 2\Omega^2$. 

Далее будет использоваться в алгоритме константу $L_0 \geq 0$, которая имеет смысл предположительной ``локальной'', константы Липшица градиента в точке $x_0$. Рассмотрим адаптивный зеркальный вариант метода подобных треугольников со стохастическим $(\delta, L)$-оракулом:
\begin{algorithm_heading}
    \label{Alg1} 
    Алгоритм \ref{Alg1} (Теоретический вариант).
\end{algorithm_heading}

\textbf{Дано:} $x_0$ -- начальная точка, $\epsilon$ -- желаемая точность решения, $\delta$, $L$ -- константы из $(\delta, L)$-оракула, $\beta$ -- доверительный уровень, $L_0$ ($L_0 \leq L$), $\sigma^2$ -- константа из определения \ref{defdeltaLstoch}.

\textbf{Алгоритм:} Возьмем 
\begin{gather*}
N := \left\lceil\frac{2\sqrt{3}\sqrt{L}R_Q}{\sqrt{\epsilon}}\right\rceil, \,
\Omega := \sqrt{6\ln{\frac{N}{\beta}}}.
\end{gather*}

\textbf{0 - шаг:}
\begin{gather*}
y_0 := x_0,\,
u_0 := x_0,\,
\alpha_0 := 0,\,
A_0 := \alpha_0.
\end{gather*}

\textbf{$\boldsymbol{k+1}$ - шаг:}

\begin{gather*} j_{k+1} := 0.\end{gather*}

\textbf{Пока не выполнится неравенство}
\begin{equation}
\begin{gathered}
\label{exitLDLST}
f_\delta(x_{k+1}) \leq f_\delta(y_{k+1}) + \langle \widetilde{\nabla}^{m_{k+1}} f_\delta(y_{k+1}), x_{k+1} - y_{k+1} \rangle\ +\\  + \frac{L_{k+1}}{2}\norm{x_{k+1} - y_{k+1}}^2 + \frac{3\sigma^2\widetilde{\Omega}}{L_{k+1}m_{k+1}} + \delta
\end{gathered}
\end{equation}

\textbf{повторять:}
\begin{gather*}
L_{k+1} := 2^{j_{k+1} - 1}L_{k}, \quad \alpha_{k+1} := \frac{1 + \sqrt{1 + 4A_kL_{k+1}}}{2L_{k+1}},\quad A_{k+1} := A_k + \alpha_{k+1},\\
y_{k+1} := \frac{\alpha_{k+1}u_k + A_k x_k}{A_{k+1}},\quad m_{k+1} := \left\lceil\frac{3\sigma^2\widetilde{\Omega}\alpha_{k+1}}{\epsilon}\right\rceil.
\end{gather*}
\begin{center}
Сгенерировать i.i.d. $\xi_j$ ($j=1,...,m_{k+1}$) и посчитать $ \widetilde{\nabla}^{m_{k+1}} f_\delta(y_{k+1})$.
\end{center}
\begin{equation*}
\begin{gathered}
\phi_{k+1}(x) \myeq V(x, u_k) + \alpha_{k+1}\left(f_\delta(y_{k+1}) + \langle \widetilde{\nabla}^{m_{k+1}} f_\delta(y_{k+1}), x - y_{k+1} \rangle\right)\\
u_{k+1} := \argmin_{x \in Q}\phi_{k+1}(x), \quad x_{k+1} := \frac{\alpha_{k+1}u_{k+1} + A_k x_k}{A_{k+1}}, \quad j_{k+1} := j_{k+1} + 1.
\end{gathered}
\end{equation*}

\textbf{Конец алгоритма.}

\begin{hypothesis}
\label{hypothesis_main}
Пусть $x_*$~--- ближайшая точка минимума к точке $x_0$ в смысле дивергенции Брэгмана, точность $\delta = \mathcal{O}(\epsilon^\frac{3}{2}).$ Тогда для предложенного алгоритма с вероятностью $1 - \mathcal{O}(\beta)$ выполнено равенство:
\begin{gather*}
f(x_N) - f(x_*) = \mathcal{O}(\epsilon).
	\end{gather*}
\end{hypothesis}


\section{Обоснование алгоритма}
\label{explanation}

Приведем обоснование данного алгоритма. Предложенный алгоритм базируется на следующих градиентных методах \cite{devolder2013exactness, gasnikov2019fast,gasnikov2017universal}. Для \textbf{адаптивного} быстрого градиентного метода с $(\delta, L)$--оракулом можно получить следующие оценки сходимости \cite{devolder2013exactness, gasnikov2019fast, gasnikov2017universal}:
\begin{equation*}
    f(x_N) - f(x_*) \leq \mathcal{O}\left(\frac{L R_Q^2}{N^2} + N \delta\right).
\end{equation*}
В то время, как для \textbf{стохастического} быстрого градиентного метода с $(\delta, L)$--оракулом можно получить оценки вида \cite{devolder2013exactness}:
\begin{equation*}
    f(x_N) - f(x_*) \leq \mathcal{O}\left(\frac{L R_Q^2}{N^2} + \frac{\sigma R_Q}{\sqrt{N}}\right).
\end{equation*}
Соответственно, ожидаемая оценка для \textbf{адаптивного стохастического} быстрого градиентного метода с $(\delta, L)$--оракулом должна иметь вид:
\begin{equation*}
    f(x_N) - f(x_*) \leq \mathcal{O}\left(\frac{L R_Q^2}{N^2} + N \delta + \frac{\sigma R_Q}{\sqrt{N}}\right).
\end{equation*}

Основная проблема по доказательству Гипотезы \ref{hypothesis_main}, которая возникает в данный момент, связана с проблемой ``выборки с отклонением'' \cite{robert2013monte}. В доказательствах методов оптимизации с стохастическим градиентом (см. \eqref{stoch_gradient}) ключевым образом используется тот факт, что стохастический градиент имеет несмещенное математическое ожижание по отношению к настоящему градиенту. В алгоритме из раздела \ref{sec:mmtDLST} это не так. Стоит обратить внимание, что наш метод является адаптивным и настраивается на локальную гладкость, делая проверку \eqref{exitLDLST}. Получается, что во внутреннем цикле, когда мы подбираем шаг метода, мы сначала генерируем mini-batch, а затем mini-batch проверяем на то, подходит ли он в неравенстве \eqref{exitLDLST}. Таким образом, мы неявно делаем процедуру ``выборки с отклонением'' \cite{robert2013monte} и соответственно возникает смещение сгенерированного стохастического градиента.

Стоит отметить, что возникали различные попытки перестроения алгоритма, но нам так и не удалось решить проблему ``выборки с отклонением''.

В алгоритме \ref{Alg1} из раздела \ref{sec:mmtDLST}, в процедуре подбора параметра $L_{k+1}$, мы каждый раз генерируем мини-батч. Возможный способ решения проблемы ``выборки с отклонением'' --- это генерировать мини-батч один раз до процедуры подбора параметра $L_{k+1}$ (см. алгоритм \ref{main_alg} из раздела \ref{comp_exp}). К сожалению, данный способ описания алгоритма приводит к другой проблеме. По аналогии с \cite{devolder2013exactness} в процессе доказательства возникает выражение вида:
\begin{equation*}
    \sum_{k = 0}^{N-1}\alpha_{k+1}\langle \widetilde{\nabla}^{m_{k+1}} f_\delta(y_{k+1})-\nabla f_\delta(y_{k+1}), x - u_k \rangle.
\end{equation*}
Основной задачей является найти некоторые хорошие оценки сверху на данную сумму. В работе \cite{devolder2013exactness} данные оценки получены и доказано, что с ``большой вероятностью'' верно следующее неравенство:
\begin{equation*}
    \sum_{k = 0}^{N-1}\alpha_{k+1}\langle \widetilde{\nabla}^{m_{k+1}} f_\delta(y_{k+1})-\nabla f_\delta(y_{k+1}), x - u_k \rangle \leq  \mathcal{O}\left(\frac{\sigma R_Q}{\sqrt{N}}\right).
\end{equation*}
Для доказательство данного неравенства использовалась неравенства типа Азума--Хефдинга (\cite{lan2012validation}, \cite{devolder2013exactness}, Лемма 7.11) и тот факт, что $\alpha_{k+1}$ являются неслучайными и условное математическое ожидание мини-батча равно настоящему градиенту. Для случая, когда мы один раз фиксируем мини-батч возникает проблема с тем, что $\alpha_{k+1}$ получаются случайные и скоррелированные с мини-батчем, поэтому стандартные неравенства типа Азума--Хефдинга ( \cite{devolder2013exactness}, Лемма 7.11) нам применить не удается.

В следующем разделе мы проводим численные эксперименты на различных задачах оптимизации. Мы показываем, что на практике практический вариант предложенного алгоритма (см. алгоритм \ref{main_alg} из раздела \ref{comp_exp}) работает лучше, чем популярные алгоритмы Adagrad \cite{duchi2011adaptive} и Adam \cite{adam}.

\newpage
\renewcommand{\algorithmicrequire}{\textbf{Input:}}
\renewcommand{\algorithmicensure}{\textbf{Output:}}

\section{Численные эксперименты}
\begin{figure}[!htb]
    \centering
\begin{subfigure}{0.4\textwidth}
  \includegraphics[width=\linewidth]{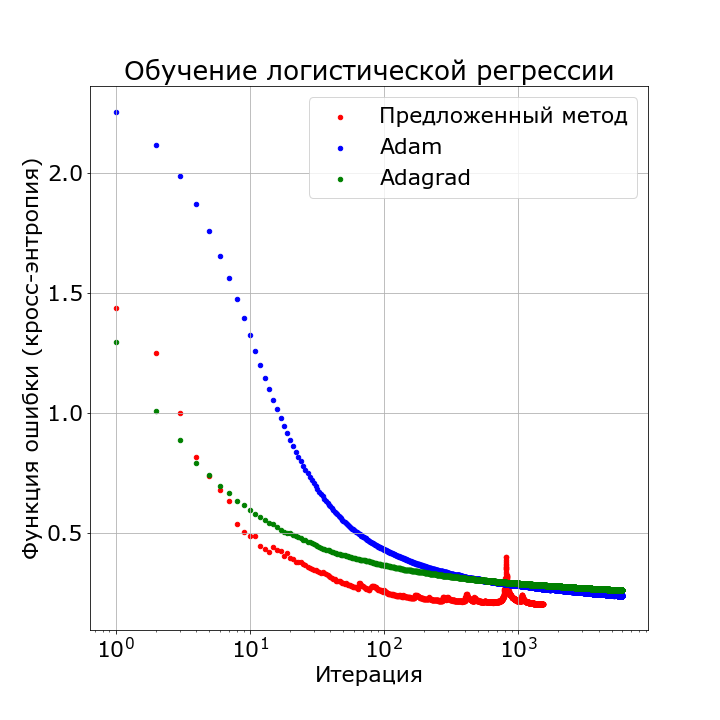}
  \caption{Функция ошибки от итерации}
  \label{fig:logobjective}
\end{subfigure}\hfil 
\begin{subfigure}{0.4\textwidth}
  \includegraphics[width=\linewidth]{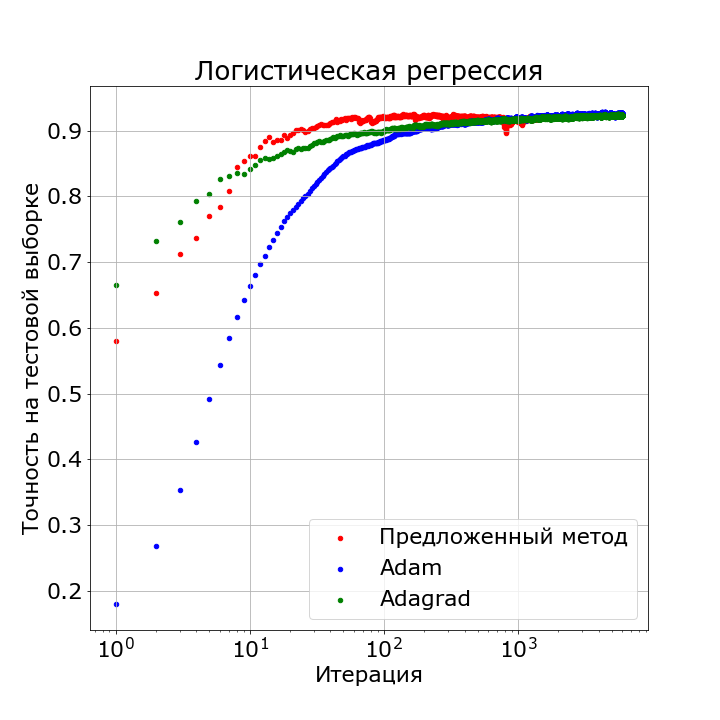}
  \caption{Тестовая точность от итерации}
  \label{fig:logtest}
\end{subfigure}

\medskip
\begin{subfigure}{0.4\textwidth}
  \includegraphics[width=\linewidth]{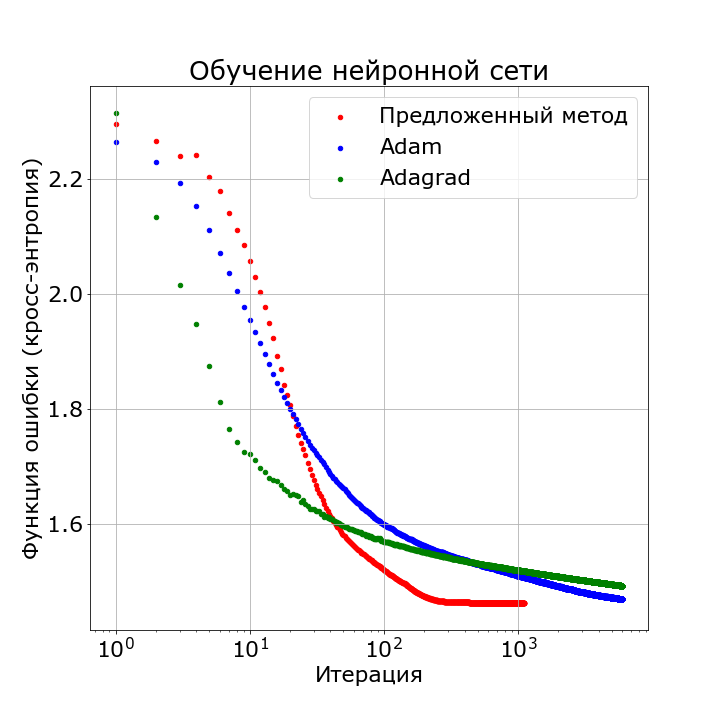}
  \caption{Функция ошибки от итерации}
  \label{fig:objective}
\end{subfigure}\hfil
\begin{subfigure}{0.4\textwidth}
  \includegraphics[width=\linewidth]{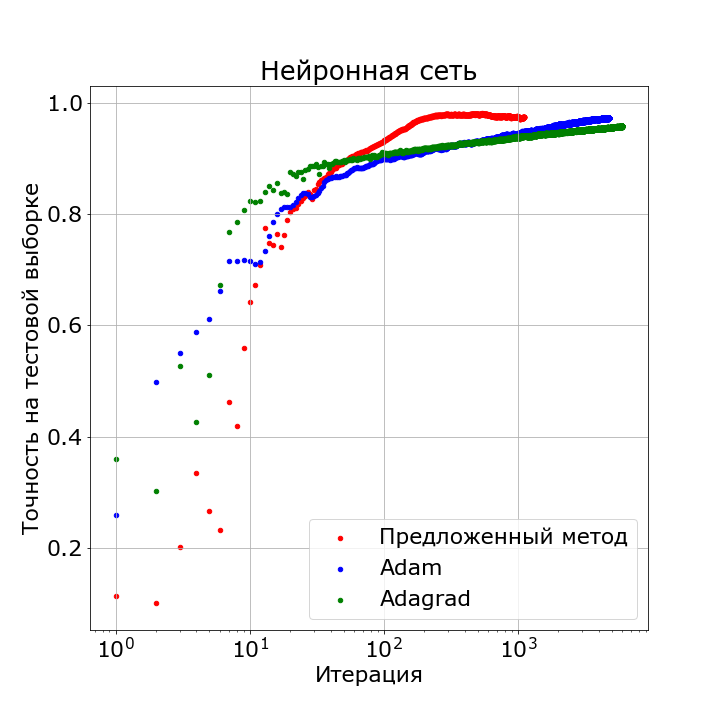}
  \caption{Тестовая точность от итерации}
  \label{fig:test}
\end{subfigure}

\medskip
\begin{subfigure}{0.4\textwidth}
  \includegraphics[width=\linewidth]{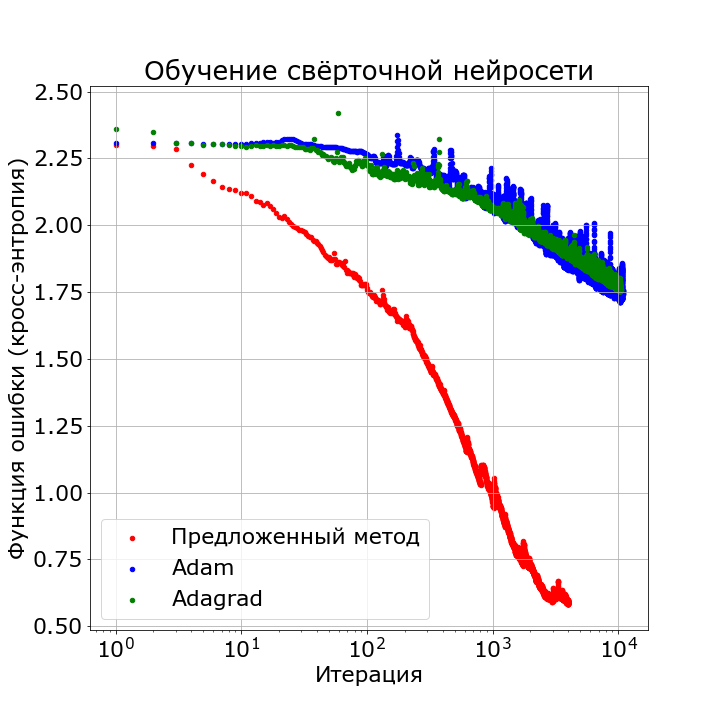}
  \caption{Функция ошибки от итерации}
  \label{fig:cnnobjective}
\end{subfigure}\hfil
\begin{subfigure}{0.4\textwidth}
  \includegraphics[width=\linewidth]{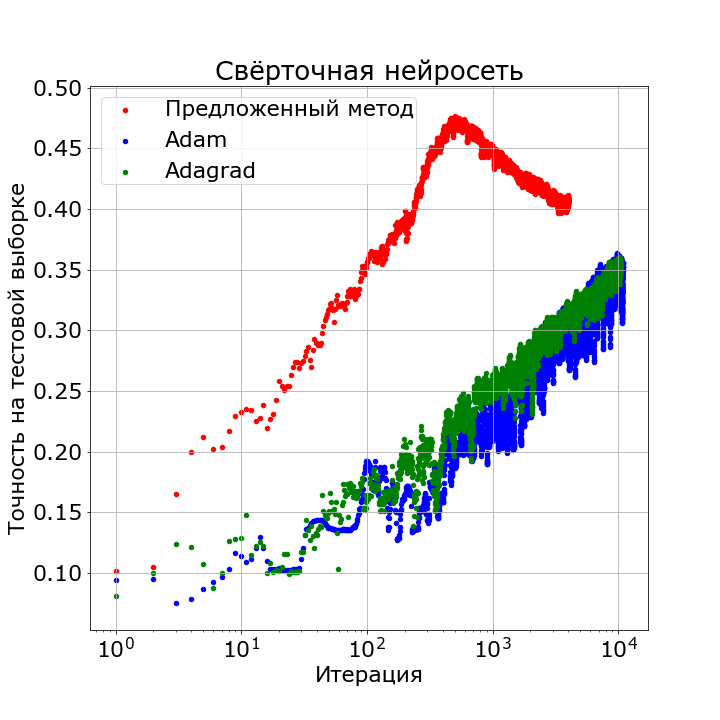}
  \caption{Тестовая точность от итерации}
  \label{fig:cnntest}
\end{subfigure}
\caption{Численные эксперименты}
\label{fig:numericexps}
\end{figure}

Начнем с описания отличия теоретического варианта (алгоритм \ref{Alg1}) и практического варианта (алгоритм \ref{main_alg}):
\begin{enumerate}
    \item На практике мы не всегда можем знать конастанту $L$ и более того константа $R_Q$ может быть не ограничена, поэтому заранее точно определить $N$ из алгоритма \ref{Alg1} мы не можем, поэтому на практике будем предполагать, что $\widetilde{\Omega} = 1$.
    \item В практическом варианте константу $m_{k+1}$ мы оцениваем, используя $L_{k}$, а не константу $L_{k+1}$. Это позволяет нам эффективнее делать подбор параметра $L_{k+1}$ за счет того, что количество $m_{k+1}$ слагаемых при подсчете mini-batch не меняется.
    \item Для многих классов оптимизируемых функций, которые включают в себя логистическую регрессию, нейронные сети и сверточные нейронные сети, сложность подсчета функции имеет тот же порядок сложности, что и подсчет градиента \cite{griewank1989automatic}. Соответственно, значение функции так же оценивается, используя технику mini-batch (см. замечание \ref{remark_stoch_fun}):
    \begin{gather*}
    f^{m_{k+1}}_\delta(y) \myeq \frac{1}{m_{k+1}}\sum_{j=1}^{m_{k+1}}\nabla f_\delta(y;\xi_j).
    \end{gather*}
\end{enumerate}
\label{comp_exp}
Теоретические рассмотрения предыдущих глав статьи приводят к практической версии алгоритма \ref{main_alg}.  Было проведено три серии экспериментов. \at2{Оптимизация всех функций происходило относительно евклидовой нормы $\norm{\cdot}_2$ c $V(x,y) = \frac{1}{2}\norm{x - y}_2^2$, множество $Q = \mathds{R}^n$.} Начальные параметры алгоритма \ref{main_alg} для всех экспериментов оставались одинаковыми: \at2{верхняя} оценка дисперсии градиента \at2{$\sigma_0^2 = 0.1$}, точность $\epsilon = 0.002$, константа $L_{0} = 1$. В экспериментах использовались следующие наборы данных:
\begin{enumerate}
    \item MNIST \cite{lecun-mnisthandwrittendigit-2010}, состоящий из изображений рукописных цифр размера $28 \times 28$ пикселей. 
    \item CIFAR \cite{cifar10}, состоящий из цветных изображений размера $3 \times 32 \times 32$ пикселей, каждое из которых принадлежит к одному из 10 классов.
\end{enumerate}
В первых двух сериях экспериментов использовался набор данных MNIST. Эксперименты заключаются в обучении логистической регрессии (линейный классификатор с выпуклой функцией потерь) и двухслойной нейронной сети с размером скрытых слоёв 1000 (нелинейный классификатор с невыпуклой функцией потерь). В случае логистической регрессии количество оптимизируемых параметров равнялось 7850, а в случае двухслойной нейронной сети уже 795010. Оптимизация проводилось при помощи предложенного алгоритма, а также при помощи наиболее распространённых на сегодняшний день адаптивных алгоритмов Adam \cite{adam} и AdaGrad \cite{duchi2011adaptive}. Для Adagrad и Adam использовались стандартные параметры: размер батча 128, $lr = 0.001$. Сравнение проводилось по скорости изменения функции потерь от итерации обучения и по изменению точности классификации на тестовой выборке. Результаты сравнения для логистической регрессии находятся на рис. \ref{fig:logobjective}, \ref{fig:logtest}, а для нейронной сети на рис. \ref{fig:objective}, \ref{fig:test}. В обоих случаях предложеный метод превосходит упомянутые выше, как по скорости сходимости функции ошибки, так и по достигнутой точности. \\
В третьей серии экспериментов небольшая свёрточная нейронная сеть обучалась предсказанию класса картинки на наборе данных CIFAR. Использовалась следующая архитектура свёрточной нейросети:
\begin{enumerate}
    \item По матрице изображения ($3 \times 32 \times 32$) проходим окном свёртки размера $3 \times 6 \times 5$. Веса свёртки являются настраиваемыми параметрами.
    \item После свёртки по полученному изображению проходит так называемый MaxPooling --- это окно размера $2 \times 2$, которое оставляет максимальный попавший в него элемент, а остальные зануляет.
    \item Далее снова проходим по результату предыдущего шага окном свёртки размера $6 \times 16 \times 5$.
    \item Далее идут три последовательных полносвязных слоя с убывающим размером: 400, 120, 84. 
\end{enumerate}
Общее число параметров данной нейросети --- 62000. В этой серии экспериментов результаты работы предложенного метода также сравнивались с Adagrad и Adam (рис. \ref{fig:cnnobjective}, \ref{fig:cnntest}). Алгоритм \ref{main_alg} также показал лучшие результаты.

\begin{algorithm_heading}
    \label{main_alg}
    Алгоритм \ref{main_alg} (Практический вариант).
\end{algorithm_heading}

\textbf{Дано:} $x_0$ -- начальная точка, константы $\epsilon > 0$, $L_0 > 0$ и $\sigma_0^2 > 0$.

\textbf{Алгоритм:}

\textbf{0 - шаг:}
\begin{gather*}
y_0 := x_0,\,
u_0 := x_0,\,
\alpha_0 := 0,\,
A_0 := \alpha_0.
\end{gather*}

\textbf{$\boldsymbol{k+1}$ - шаг:}

\begin{gather*} \tilde{\alpha}_{k+1} = \frac{1 + \sqrt{1 + 4A_{k} L_{k}}}{2L_{k}}, \quad m_{k+1} := \left\lceil\frac{3\sigma_0^2\tilde{\alpha}_{k+1}}{\epsilon}\right\rceil, \quad j_{k+1} := 0,\end{gather*}
\begin{center}
сгенерировать i.i.d. $\xi_j$ ($j=1,...,m_{k+1}$).
\end{center}

\textbf{Пока не выполнится неравенство}
\begin{equation}
\begin{gathered}
f^{m_{k+1}}_\delta(x_{k+1}) \leq f^{m_{k+1}}_\delta(y_{k+1}) + \langle \widetilde{\nabla}^{m_{k+1}} f_\delta(y_{k+1}), x_{k+1} - y_{k+1} \rangle\ +\\  + \frac{L_{k+1}}{2}\norm{x_{k+1} - y_{k+1}}^2 + \frac{\epsilon}{L_{k+1}\alpha_{k+1}}
\end{gathered}
\end{equation}

\textbf{повторять:}
\begin{gather*}
L_{k+1} := 2^{j_{k+1} - 1}L_{k}, \quad \alpha_{k+1} := \frac{1 + \sqrt{1 + 4A_kL_{k+1}}}{2L_{k+1}},\quad A_{k+1} := A_k + \alpha_{k+1},\\
y_{k+1} := \frac{\alpha_{k+1}u_k + A_k x_k}{A_{k+1}}.
\end{gather*}
\begin{center}
Посчитать $ \widetilde{\nabla}^{m_{k+1}} f_\delta(y_{k+1})$ и $f^{m_{k+1}}_\delta(y_{k+1})$.
\end{center}
\begin{equation*}
\begin{gathered}
\phi_{k+1}(x) \myeq V(x, u_k) + \alpha_{k+1}\left(f^{m_{k+1}}_\delta(y_{k+1}) + \langle \widetilde{\nabla}^{m_{k+1}} f_\delta(y_{k+1}), x - y_{k+1} \rangle\right)\\
u_{k+1} := \argmin_{x \in Q}\phi_{k+1}(x), \quad x_{k+1} := \frac{\alpha_{k+1}u_{k+1} + A_k x_k}{A_{k+1}}, \quad j_{k+1} := j_{k+1} + 1.
\end{gathered}
\end{equation*}
\begin{center}
Посчитать $f^{m_{k+1}}_\delta(x_{k+1})$.
\end{center}

\textbf{Конец алгоритма.}

\newpage

\printbibliography

\end{document}